\documentclass [a4paper,10pt]{article}
\usepackage [utf8] {inputenc}
\usepackage {amsmath}
\usepackage {amsthm}
\usepackage {amssymb}
\usepackage {amsfonts}
\newtheorem {lemma}{Lemma}
\newtheorem* {theorem}{Theorem}

\title {Primes in arithmetic progressions on average I}
\author {Tomos Parry}
\date {}
\begin {document}
\maketitle
\begin {abstract}
Let
\begin {eqnarray*}
E_x(q,a)&:=&\sum _{p\leq x\atop {p\equiv a(q)}}\log p-\frac {x}{\phi (q)}\hspace {10mm}M(Q):=\sum _{q\leq Q}\phi (q)\sideset {}{'}\sum _{a=1}^qE_x(q,a)^3.
\end {eqnarray*}
We show that if $Q$ is close to $x$ then $M(Q)\ll Q^3(x/Q)^{7/5}$, thereby showing that sign-changes in the average gives power-saving cancellation past the expected $\sqrt {x/q}$. This therefore for the first time provides unconditional evidence for a deep conjecture of Montgomery and Soundararajan.
\end {abstract}
\begin {center}
\section {Introduction}
\end {center}
Let
\[ E_x(q,a):=\sum _{p\leq x\atop {p\equiv a(q)}}\log p-\frac {x}{\phi (q)}\]
be the error term when counting primes in arithmetic progressions,and which (ignoring $\log $ factors) is famously suspected to be around $\sqrt {x/q}$ in size. When averaging, however, additional cancellation may occur, as was demonstrated by Hooley \cite {hooley} in the form
\begin {eqnarray*}
M(Q)&:=&\sum _{q\leq Q}\phi (q)\sideset {}{'}\sum _{a=1}^qE_x(q,a)^3
\\ &\ll &Q^3\left (\frac {x}{Q}\right )^{3/2}e^{-c\sqrt {\log x}}+\frac {x^3}{(\log x)^A}.
\end {eqnarray*}
Hooley remarks further ``\emph {The assumption of the extended Riemann hypothesis for Dirichlet's $L$-functions enables one to obtain a greatly improved form of our result in which a sharper version of the above bound is valid for smaller values of $Q$, thus shedding further light on the distribution of $E_x(q,a)$ for smaller values of $q$. But the investigation of this development must await a later paper in this series.}". It seems Hooley never did manage to return to this development but following his proof we see that he ultimately reduces $M(Q)$ to an integral looking something like
\[ Q^3\int _{(-1/2)}\frac {\zeta (s)\zeta (s+1)h(s)(x/Q)^{s+2}ds}{\zeta (2s+2)s(s+1)(s+2)}\]
with $h(s)$ holomorphic and bounded for $\sigma >-1$, so we'd guess he meant that under RH the above integral may be moved to the left (as well as stronger bounds then being available for the circle method part of his argument) to imply something like
\[ M(Q)\ll Q^3\left (\frac {x}{Q}\right )^{5/4}+x^{5/2}\]
although, in respect to the wealth of his contributions, we should also allow for the possibility that he meant something more.
\\
\\ In any case, as far as we know, there is no trace of anything in the literature giving additional cancellation of $E_x(q,a)$ past the expected $\sqrt {x/q}$ heuristic further than Hooley's above result, and certainly no power-saving without RH. In this paper, we reveal an additional cancellation within the dispersion method which allows us to (unconditionally) prove
\\
\begin {theorem} 
For $Q\leq x/(\log x)^2$ and any $A>0$
\[ M(Q)\ll _AQ^3\left (\frac {x}{Q}\right )^{7/5}+\frac {x^3}{(\log x)^A}.\]
\end {theorem}
\hspace {1mm}
\\ Note that a deep conjecture of Montgomery proposes
\[ \int _1^X\left (\sum _{x<n\leq x+H}\Lambda (n)-H\right )^2dx\sim XH\log (X/H)\hspace {10mm}X^\epsilon \leq H\leq X^{1-\epsilon }\]
to which Montgomery and Soundararajan \cite {montsound} have shown, on the assumption of some heavy conjectures, 
\[ \int _1^X\left (\sum _{x<n\leq x+H}\Lambda (n)-H\right )^kdx\sim XH^{k/2}\log (X/H)^{k/2}\left (\mu _k+\mathcal O\left (H^{-c/k}\right )\right )\hspace {10mm}X^\epsilon \leq H\leq X^{1/k}\]
(for some $c>0$, and where the $\mu _k$ are the moments of the standard Gaussian). The $q$-analogue of this (and $q$-analogues are typically certainly no easier) for $k=3$ would read
\[ \sideset {}{'}\sum _{a=1}^qE_x(q,a)^3\ll q(x/q)^{3/2-c}\]
and in fact Leung \cite {leung} has proven a result of this kind, again highly conditional. Our result therefore, as far as we know, is the first that detects the conjectured cancellation unconditionally.
\\
\\ Let's give a few words as to the proof. Opening up $M(Q)$ we see it's
\[ \sum _{q\leq Q}\phi (q)\sum _{p_1,p_2,p_3\leq x\atop {p_1\equiv p_2\equiv p_3(q)}}\log \mathbf p-3x\sum _{q\leq Q}\sideset {}{'}\sum _{a=1}^qE_x(q,a)^2-3x^2E_x(1,0)\sum _{q\leq Q}\frac {1}{\phi (q)}-x^3\sum _{q\leq Q}\frac {1}{\phi (q)^2}\]
where the first term is ``the difficult part" and where let's call the second term ``the variance part". Hooley finds a precise asymptotic for the difficult part and compares it to the main terms of the remaining parts, including the variance part. We prove our theorem following this argument, but instead of bounding the integral above, we further compare it to the \emph {error} term in the variance part, and show they completely negate each other. This requires, and what is most likely the most important contribution of our proof, keeping better track of the main terms than in \cite {hooley}, which in (117) contains a large term impossible to control without RH. This laying-out of the main terms in a more transparent manner so as to isolate the correct main term is essential if we wish to pass to the study of the fourth moment with a dispersion method, and it is indeed that subject which is the main motivation for the present work. We return to this in a sequel to this article, but let's just briefly indicate now what we mean: Opening up the fourth moment as was done above for the third moment, and using (e.g.) a weight $\phi (q)^2$, we get a ``difficult part" plus parts
\[ -4x\underbrace {\sum _{q\leq Q}\phi (q)\sideset {}{'}\sum _{a=1}^qE_x(q,a)^3}_{=M(Q)}+6x^2\sum _{q\leq Q}\sideset {}{'}\sum _{a=1}^qE_x(q,a)^2+...\]
so that if we want a bound for the fourth moment of size $\ll Q^4(x/Q)^{2+\delta }$ then at the very least we need a bound $M(Q)\ll Q^3(x/Q)^{1+\delta }$.
\\ 
\\ As a final but ultimately just vague comment, we mention that error terms in integrals negating each other when otherwise needing to call on RH is something we believe should maybe not be glossed over.
\\
\\ We present our proof in two sections. In the first we prove our theorem subject to Lemma \ref {j} - here we are just re-writing the first half of \cite {hooley}. In the second we state and prove Lemma \ref {j}, where we find the new cancellation. In particular this lemma gives us a contour integral with integrand containing a factor $1/s^2(s+1)$ instead of the $1/s(s+1)(s+2)$ in \cite {hooley}. This seemingly completely uninteresting observation proves to be just the opposite.
\\
\\ We particularly thank Professor Ahmet Güloğlu for many helpful discussions in writing this paper. 
\begin {center}
\section {Proof of theorem, subject to Lemma \ref {j}}
\end {center}
We retrace initially Hooley's route exactly, however we give all the details since there are a few inaccuracies\footnote {Essentially \S 4 needs to be re-written, the ultimate problem being that smoothness of sums involved is overlooked. (22) to (23) misses a large term, causing a sharp cut-off, and the proof of Lemma 2 is given as a reference to another lemma but where the sum is smoother than in Lemma 2. This is all fixable and of course doesn't detract from the main point, but it probably does mean that we'd be wise to write things out again to clear up any possible blurriness.} in \cite {hooley}. Notation: By $\int _{(c)}$ we mean $\frac {1}{2\pi i}\int _{c+1/100-i\infty }^{c+{1/100}+i\infty }$.
\\
\\ Let $A>0$ and $P:=x/(\log x)^{A+100}$. By (Gallagher's form of) the Barban-Davenport-Halberstam theorem
\[ \sum _{q\leq P}\sideset {}{'}\sum _{a=1}^qE_x(q,a)^2\ll 
\frac {x^2}{(\log x)^{A}}\]
so, with the prime number theorem,
\begin {eqnarray}\label {dechra}
M(Q)&=&\sum _{P<q\leq Q}\phi (q)\sideset {}{'}\sum _{a=1}^qE_x(q,a)^3+\mathcal O\left (\frac {x^3}{(\log x)^A}\right )\notag 
\\ &=&\sum _{P<q\leq Q}\phi (q)\sum _{\mathbf p\leq x\atop {p_1\equiv p_2\equiv \tilde p_3(q)}}\log \mathbf p-3x\sum _{P<q\leq Q}\sideset {}{'}\sum _{a=1}^qE_x(q,a)^2\notag 
\\ &&-\hspace {4mm}3x^2E_x(1,0)\sum _{P<q\leq Q}\frac {1}{\phi (q)}-x^3\sum _{P<q\leq Q}\frac {1}{\phi (q)}+\mathcal O\left (\frac {x^3}{(\log x)^A}\right )\notag 
\\ &=:&S_1(Q)-3xV(Q)+\mathcal O\left (\frac {x^3}{(\log x)^A}\right )-x^3U(Q)+\mathcal O\left (\frac {x^3}{(\log x)^A}\right )
\end {eqnarray}
and by (2), (5), (6), (8), (9), (10) and (11) of \cite {hooley1}
\begin {eqnarray}\label {variance}
V(Q)&=&Qx\log x-Qx-x^2U(Q)+x^2\Big (J(x/P)-J(x/Q)\Big )+\mathcal O\left (\frac {x^2}{(\log x)^A}\right )\notag 
\\ \text {where }\hspace {5mm}J(X)&:=&Z\log X+C_1+\frac {\log X+C_2}{X}+2\int _{(-5/4)}\frac {\mathcal F_1(s+1)(x/Q)^{s}ds}{s(s+1)(s+2)}\notag 
\\ \mathcal F_1(s)&:=&\sum _{n=1}^\infty \frac {n/\phi (n)}{n^s}\notag 
\\ Z&:=&\frac {\zeta (2)\zeta (3)}{\zeta (6)}\hspace {10mm}C_2:=\frac {\zeta '(0)}{\zeta (0)}+\gamma _0+\sum _p\frac {\log p}{p(p-1)};
\end {eqnarray}
here we get the value of $C_2$ from (10) of \cite {hooley}. From Lemma 1 of \cite {hooley1} 
\begin {equation}\label {u}
U(Q)=Z\log (Q/P)+\mathcal O\left (\frac {\log P}{P}\right )
\end {equation}
and of course $J(x/P)\ll P/x$ so \eqref {variance} becomes
\begin {eqnarray}\label {variance11}
V(Q)&=&Qx\log Q-(C_2+1)Qx+2Q^2\underbrace {\int _{(-1/4)}\frac {\mathcal F_1(s)(x/Q)^{s+1}ds}{(s-1)s(s+1)}}_{=:I(x/Q)}+\mathcal O\left (Px\log x\right )
\end {eqnarray}
and with \eqref {u} then \eqref {dechra} becomes, writing $X=x/Q$,
\begin {eqnarray}\label {variance22}
M(Q)=S_1(Q)-3Q^3\Big (X^2\log Q-X^2(C_2+1)+2XI(X)\Big )-3Zx^3\log (Q/P)+\mathcal O\left (\frac {x^3}{(\log x)^A}\right ).
\end {eqnarray}
 Let
\[ J_1(Q):=\sum _{Q<q\leq x}\frac {\phi (q)}{q}\sum _{\mathbf p\leq x\atop {p_1\equiv p_2\equiv p_3(q)}}\log \mathbf p\]
so that
\begin {eqnarray}\label {parts}
S_1(Q)&=&Q\Big (J_1(P)-J_1(Q)\Big )-\int _P^Q\Big (J_1(P)-J_1(t)\Big )dt.
\end {eqnarray}
Then, using $\phi \star 1=\text { id }$ and where $f(x)$ is some quantity independent of $Q$, 
\begin {eqnarray}\label {jj}
J_1(Q)&=&
\sum _{d\leq x}\frac {\mu (d)}{d}\Bigg \{ 6\sum _{Q<q\leq x\atop {d|q}}\sum _{p_1<p_2<p_3\leq x\atop {p_1\equiv p_2\equiv p_3(q)}}\log \mathbf p+3\sum _{Q<q\leq x\atop {d|q}}\sum _{p<p'\leq x\atop {p\equiv p'(q)}}(\log p)^2\log p'\notag 
\\ &&+\hspace {4mm}3\sum _{Q<q\leq x\atop {d|q}}\sum _{p<p'\leq x\atop {p\equiv p'(q)}}\log p(\log p')^2\Bigg \} +\sum _{Q<q\leq x}\frac {\phi (q)}{q}\sum _{p\leq x}(\log p)^3
\notag 
\\ &=:&6J_4(Q)+3\Big (J_5(Q)+J_6(Q)\Big )+f(x)+\mathcal O\left (xQ(\log x)^3\right )
\end {eqnarray}
so we need to calculate the $J_i$'s - first we do $J_5$ and $J_6$. 
Let 
\begin {eqnarray*}
\psi _1(l)&:=&\prod _{p|l}\left (1+\frac {1}{p-1-1/p}\right )\hspace {5mm}\mathcal F(s):=\sum _{l=1}^\infty \frac {\psi _1(l)}{l^s}\hspace {5mm}\text {so }\hspace {5mm}\mathcal F(s)=\zeta (s)\zeta (s+1)\underbrace {\mathcal F^*(s)}_{\text {$\ll 1$ for $\sigma >-1/2$}}=:\zeta (s)h(s)
\\ \alpha _1&=&h(1)\hspace {10mm}\beta _1=h(1)(\gamma _0-1)+h'(1)\hspace {10mm}
\\ \alpha _2&=&\frac {h(1)}{2}\hspace {10mm}\beta _2=h(1)\left (\frac {\gamma _0}{2}-\frac {3}{4}\right )+\frac {h'(1)}{2}
\end {eqnarray*}
so that, from the usual Perron type formulas and moving the integrals a bit past $\sigma =-1$ with the simplest of bounds for $\zeta (s)$,
\begin {eqnarray*}
\sum _{l<X}(X-l)\frac {\psi _1(l)}{l}&=&\int _{(0)}\frac {\mathcal F(s+1)X^{s+1}ds}{s(s+1)}
=X\left (\alpha _1\log X+\beta _1\right )+\mathcal O(1)\hspace {10mm}
\\ \sum _{l<X}(X-l)^2\frac {\psi _1(l)}{l}&=&2\int _{(0)}\frac {\mathcal F(s+1)X^{s+2}ds}{s(s+1)(s+2)}=2X^2\left (\alpha _2\log X+\beta _2\right )+\mathcal O(X)
\\ \sum _{l<X}(X-l)\psi _1(l)&=&\int _{(0)}\frac {\mathcal F(s+1)X^{s+2}ds}{(s+1)(s+2)}
=c_1X^2+\mathcal O(X)
\\ \sum _{l<X}\log (X/l)l\psi _1(l)&=&\int _{(0)}\frac {\mathcal F(s+1)X^{s+2}ds}{(s+2)^2}=c_2X^2+\mathcal O(X)\hspace {10mm}\text {for some constants $c_1,c_2$}.
\end {eqnarray*}
Let
\begin {eqnarray*}
\mathcal M_1(l)&:=&Q^2(X-l)\left (\gamma X+\delta l\right )-\frac {(lQ)^2\log (X/l)}{2}\hspace {10mm}X:=\frac {x}{Q}\hspace {6mm}\gamma :=\frac {\log x}{2}-\frac {1}{4}\hspace {6mm}\delta =\frac {3}{4}-\frac {\log x}{2}
\\ \mathcal M_2(l)&:=&\int _0^{x-Ql}\log t\hspace {1mm}(x-Ql-t)dt
\end {eqnarray*}
and note that
\begin {eqnarray}\label {mu}
\sum _{d\leq x}\frac {\mu (d)}{d}\frac {1}{\phi (dl)}&=&
\frac {C_3\psi _1(l)}{l}+\mathcal O\left (\frac {\log x}{x}\right )\hspace {7mm}\text {where}\hspace {7mm}C_3=\prod _{p}\left (1-\frac {1}{p(p-1)}\right )
\end {eqnarray}
so the integral representations above give (again for some $f(x)$ independent of $Q$)
\begin {eqnarray}\label {j5j6}
\sum _{l\leq X}\frac {\mathcal M_1(l)\psi _1(l)}{l}&=&\gamma \alpha _1x^2\log X+f(x)+\mathcal O\left (Q^2X\right )\notag 
\\ \sum _{l\leq X}\frac {\mathcal M_2(l)\psi _1(l)}{l}&=&\frac {Q^2\log Q}{2}\sum _{l\leq X}(X-l)^2\frac {\psi _1(l)}{l}+Q^2\int _0^X\log (X-t)\left (\sum _{l<t}(t-l)\frac {\psi _1(l)}{l}\right )dt\notag 
\\ &=&Q^2X^2\log Q\Big (\alpha _2\log X+\beta _2\Big )
\notag 
\\ &&\hspace {4mm}+\hspace {4mm}Q^2\int _0^Xt\log (X-t)\Big \{ \alpha _1\log t+\beta _1\Big \} dt\mathcal O\left (Q^2X\log X\right )\notag 
\\ &=&x^2\log Q\Big (\alpha _2\log X+\beta _2\Big )\notag   
\\ &&\hspace {4mm}\hspace {4mm}x^2\left (\frac {\alpha _1(\log X)^2}{2}+\left (\frac {\beta _1}{2}-\alpha _1\right )\log X\right )+f(x)+\mathcal O(Q^2X\log X);
\end {eqnarray}
if this last integral proves tricky to check we can always refer to the proof of Lemma 2 in \cite {hooley} and read from line -6 to the end. Recall $Q\geq x/(\log x)^{A+100}$ and put $L=(\log x)^{A+10}x/Q$. We have from the Siegel-Walfisz theorem
\begin {eqnarray*}
\sum _{l<x/Q}\sum _{p,p'\leq x\atop {p\equiv p'(dl)\atop {Q<\frac {p'-p}{l}}}}\log p(\log p')^2&=&\sum _{l<x/Q}\sum _{Ql<p'\leq x}(\log p')^2\left (\frac {p'-Ql}{\phi (dl)}+\mathcal O\left (\frac {x}{L}\right )\right )+\mathcal O\left (\frac {x^{2+\epsilon }}{Q}\right )
\\ &=&\sum _{l<x/Q}\frac {1}{\phi (dl)}\Bigg \{ \int _{Ql}^x(t-Ql)\log t\hspace {1mm}dt+\mathcal O\left (\frac {x^2}{(\log x)^{A+2}}\right )\Bigg \} +\mathcal O\left (\frac {x^2}{(\log x)^{A+1}}\right )
\\ &=&\sum _{l<x/Q}\frac {\mathcal M_1(l)}{\phi (dl)}+\mathcal O\left (\frac {x^2}{(\log x)^{A+1}}\right )
\end {eqnarray*}
and similarly
\begin {eqnarray*}
\sum _{l<x/Q}\sum _{p,p'\leq x\atop {p\equiv p'(dl)\atop {Q<\frac {p'-p}{l}}}}(\log p)^2\log p'&=&\sum _{l<x/Q}\sum _{p\leq x-Ql}(\log p)^2\left (\frac {x-Ql-p}{\phi (dl)}+\mathcal O\left (\frac {x}{L}\right )\right )+\mathcal O\left (\frac {x^{2+\epsilon }}{Q}\right )
\\ &=&\sum _{l<x/Q}\frac {1}{\phi (dl)}\left (\mathcal M_2(l)+\mathcal O\left (\frac {x^2}{(\log x)^{A+2}}\right )\right )+\mathcal O\left (\frac {x^2}{(\log x)^{A+1}}\right )
\end {eqnarray*}
so, re-writing the congruence condition, $J_5(Q)+J_6(Q)$ is
\begin {eqnarray*}
\sum _{d\leq x}\frac {\mu (d)}{d}\left \{ \sum _{l<x/Q}\frac {1}{\phi (dl)}\Bigg (\mathcal M_1(l)+\mathcal M_2(l)\Bigg )+\mathcal O\left (\frac {x^2}{(\log x)^{A+1}}\right )\right \} 
\end {eqnarray*}
and we conclude from \eqref {mu} and \eqref {j5j6}
\begin {eqnarray}\label {bb}
J_5(Q)+J_6(Q)&=&
C_3x^2\left (\frac {\alpha _1(\log X)^2}{2}+\left (\gamma \alpha _1+\frac {\beta _1}{2}-\alpha _1\right )\log X+\alpha _2\log Q\log X+\beta _2\log Q\right )\notag 
\\ &&+\hspace {4mm}f(x)+\mathcal O\left (Q^2X\log X+\frac {x^2}{(\log x)^A}\right )\notag 
\\ &=&C_3x^2\log Q\left (\left (\alpha _2-\frac {3\alpha _1}{2}\right )\log x+\left (\frac {5\alpha _1}{4}-\frac {\beta _1}{2}+\beta _2\right )\right )+f(x)+\mathcal O\left (Qx\log x+\frac {x^2}{(\log x)^A}\right )\notag 
\\ &=&x^2\log Q\Big (-\log x+1\Big )+f(x)+\mathcal O\left (Qx\log x+\frac {x^2}{(\log x)^A}\right )\hspace {10mm}(\text {because }C_3h(1)=1)\notag 
\\ &=:&\mathcal B(Q)+f(x)+\mathcal O\left (Qx\log x+\frac {x^2}{(\log x)^A}\right ).
\end {eqnarray}
Now to $J_4$. Whenever the letters $l,l'$ are in context write $\tilde l=l+l'$. Let
\begin {eqnarray}\label {defn}
r&:=&\frac {1}{p-2-\frac {1}{p}}\left \{ \begin {array}{ll}1&\text { if }p>2\\ 0&\text { if }p=2\end {array}\right \} \hspace {10mm}
\Gamma _\Delta (l):=\prod _{p|l\atop {p\nmid \Delta }}\left (1+r\right )\hspace {10mm}\notag 
\\ I(p)&:=&\left \{ \begin {array}{ll}-r(1+3r+r^2)&\text { if $p$ is odd}\\ 2&\text { if $p=2$ }\end {array}\right \} \hspace {10mm}I(\Delta ):=\prod _{p|\Delta }\left \{ \begin {array}{ll}I(p)&\text { if $\Delta $ is squarefree}\\ 0&\text { if not }\end {array}\right \} \notag 
\\ \mathcal J(X)&:=&\sum _{\Delta (l+l')\leq X}\frac {I(\Delta )(X-\Delta \tilde l)^2\Gamma _\Delta (l)\Gamma _\Delta (l')\Gamma _\Delta (\tilde l)}{\Delta \tilde l}\hspace {10mm}\notag 
\\ \mathcal J^*(X)&:=&\sum _{\Delta (l+l')\leq X}I(\Delta )(X-\Delta \tilde l)^2\Delta \tilde l\Gamma _\Delta (l)\Gamma _\Delta (l')\Gamma _\Delta (\tilde l)\notag 
\\ C_5&:=&\prod _{p>2}\left (1-\frac {1}{(p-1)^2}\right )\hspace {10mm}C_6:=\prod _{p>2}\left (1-\frac {1}{p(p-2)}\right )\hspace {10mm}C_8:=\prod _{p>2}\left (1+\frac {2r}{p}\right )\notag 
\\ \mathcal U(s)&:=&\prod _{p}\left (1+\frac {Y}{p-1}\right )\hspace {57mm}\sigma >0\hspace {12mm}Y:=p^{-s}\notag 
\\ &=&\zeta (s+1)\prod _{p}\left (1+\frac {Y}{p-1}\right )(1-Y/p)\hspace {29mm}\sigma >-1/2\notag 
\\ a&:=&\frac {2\Gamma (3)}{\Gamma (6)}=\frac {1}{30}\hspace {10mm}b=\frac {2\Gamma (3)\zeta (0)}{\Gamma (5)}=-\frac {1}{12}\hspace {10mm}Z:=\frac {\zeta '(0)}{\zeta (0)}+\gamma _0-\frac {13}{12}\notag 
\\ \alpha &:=&a\hspace {0.5mm}\mathcal U(1)\hspace {10mm}\beta :=b\hspace {0.5mm}\mathcal U(0)\hspace {10mm}\gamma :=\beta \left (Z+\sum _p\frac {\log p}{p(p-1)}\right )+\frac {\zeta (0)}{12}\notag 
\\ M(X)&:=&\alpha X^5+\beta X^4\log X+\gamma X^4\hspace {10mm}E(X):=2\int _{(-1/2)}\frac {\zeta (s)\mathcal U(s)X^{s+4}ds}{s(s+3)(s+4)}
\end {eqnarray}
and denote by $\Gamma (s)$ the Gamma function. Hooley deals with $J_4$ using the circle method, this original idea being the main point of the work. There are no inaccuracies here as far as we can tell, so we simply report that 
\[ J_4(X)=\frac {C_5C_6Q^2}{2}\mathcal J(X)+\mathcal O\left (\frac {x^2}{(\log x)^A}\right )\]
which (referring to \cite {hooley}) we can read off from (13), (18), (54), (59), (60), (63), (64) and the sentence following (64), and we get from \eqref {jj} and \eqref {bb}
\begin {eqnarray}\label {j1}
J_1(Q)&=&3C_5C_6Q^2\mathcal J(X)+3\mathcal B(Q)+f(x)+\mathcal O\left (xQ(\log x)^2+\frac {x^2}{(\log x)^A}\right ).
\end {eqnarray}
for some $f(x)$ independent of $Q$. We now need
\\
\begin {lemma}\label {j}
Assume the notation \eqref {defn}. For $\delta =1/10$
\begin {eqnarray*}
\mathcal J^*(X)&=&C_8\Big (M(X)+E(X)\Big )+\mathcal O\left (X^{7/2-\delta }\right ).
\end {eqnarray*}
\end {lemma}
This should be compared with (117) of \cite {hooley}, where the main term can't be controlled without RH and where the error term already contains the usual ``zero-free region error". 
\\
\\ We prove this lemma in the next section. Assuming it for the moment, write
\begin {eqnarray*}
\mathcal C_f(X)&=&\frac {f(X)}{X^2}-6X\int _0^{X}\frac {f(t)dt}{t^4}+12X^2\int _{0}^{X}\frac {f(t)dt}{t^5}
\end {eqnarray*}
so that, for some $f(x)$ independent of $Q$,
\begin {eqnarray*}
\mathcal C_{M}(X)&=&10\alpha X^3+6\beta X^2(\log X)^2+(12\gamma -5\beta )X^2\log X+X^2(6\beta -5\gamma )
\\ &=&\frac {\mathcal U(1)X^3}{3}-\frac {\mathcal U(0)X^2(\log X)^2}{2}+\mathcal U(0)\left (1\underbrace {-\frac {\zeta '(0)}{\zeta (0)}-\gamma _0-\sum _p\frac {\log p}{p(p-1)}}_{=-C_2\text { as in \eqref {variance}}}\right )X^2\log X
\\ &=:&AX^3+BX^2(\log X)^2+CX^2\log X
\\ &=&AX^3-X^2\Big (\log Q(2B\log x+C)-B(\log Q)^2+f(x)\Big )
\\ &=:&M^*(Q)
\end {eqnarray*}
and 
\begin {eqnarray*}
\mathcal C_{E}(X)&=&\int _{(-1/2)}\frac {\zeta (s)\mathcal U(s)}{s(s+3)(s+4)}\mathcal C_{x\mapsto x^s}(X)ds
=\int _{(-1/2)}\frac {\zeta (s)\mathcal U(s)X^{s+2}}{s^2(s+1)}ds
=:E^*(X)
\end {eqnarray*}
to conclude
\begin {eqnarray*}
P^2\mathcal C_{\mathcal J^*}(x/P)-Q^2\mathcal C_{\mathcal J^*}(x/Q)
&=&C_8\Big (P^2M^*(x/P)-Q^2M^*(x/Q)+P^2E^*(x/P)-Q^2E^*(x/Q)\Big )+\mathcal O\left (Q^{2}\left (\frac {x}{Q}\right )^{7/5}\right )\notag 
\\ &=:&C_8\Big (\mathcal A(Q)+\mathcal E(Q)\Big )+\mathcal O\left (Q^{2}\left (\frac {x}{Q}\right )^{7/5}\right ).
\end {eqnarray*}
From Lemma 3 of \cite {hooley} and \eqref {j1} this means
\begin {eqnarray*}
J_1(P)-J_1(Q)&=&
3\Big (\mathcal A(Q)+\mathcal B(Q)+\mathcal E(Q)\Big )+\mathcal O\left (Q^{2}\left (\frac {x}{Q}\right )^{7/5}+\frac {x^2}{(\log x)^A}\right )\hspace {10mm}(C_5C_6C_8=1)
\\ &=&3A\left (\frac {x^3}{P}-\frac {x^3}{Q}\right )+3x^2\Bigg (\Big ((2B+1)\log x+C-1\Big )\log (P/Q)+B(\log P)^2-B(\log Q)^2\Bigg )
\\ &&\hspace {0mm}+\hspace {4mm}3\mathcal E(Q)+\mathcal O\left (Q^{2}\left (\frac {x}{Q}\right )^{7/5}\right )\hspace {10mm}\text {(from \eqref {bb}; note also $\mathcal U(0)=1$)}
\\ &=&3A\left (\frac {x^3}{P}-\frac {x^3}{Q}\right )+3x^2\Big (C_2\log (P/Q)+B(\log P)^2-B(\log Q)^2\Big )+3\mathcal E(Q)+\mathcal O\left (Q^{2}\left (\frac {x}{Q}\right )^{7/5}\right )
\\ &=&3x^2M^{**}(Q)+3\mathcal E(Q)+\mathcal O\left (Q^{2}\left (\frac {x}{Q}\right )^{7/5}\right )
\end {eqnarray*}
so by \eqref {parts} and the above definition for $\mathcal E(Q)$
\begin {eqnarray*}
S_1(Q)&=&3x^2\Bigg \{ QM^{**}(Q)-\int _P^QM^{**}(t)dt+2\int _{(-1/2)}\frac {\zeta (s)\mathcal U(s)x^{s}}{s^2(s+1)}\underbrace {\left (Q(P^{-s}-Q^{-s})-\int _P^Q(P^{-s}-t^{-s}\right )}_{=Q^{1-s}\frac {s}{s-1}+\mathcal O(xP^{-s})}\Bigg \} .
\end {eqnarray*}
The first term in curly brackets is
\begin {eqnarray*}
&=&\int _P^Qt\left (\frac {Ax}{t^2}-\frac {S}{t}-\frac {2B\log t}{t}\right )dt
\\ &=&Ax\log (Q/P)-Q\Big (2B\log Q-2B+C_2\Big )+\mathcal O(P\log P)
\\ &=&Ax\log (Q/P)+Q\left (\log Q-1-C_2\right )
\end {eqnarray*}
and from \eqref {variance11} the second term gives a total contribution
\[ 6xQ^2I(x/Q)+\mathcal O\left (\frac {x^2}{(\log x)^A}\right )\]
and our theorem now follows from \eqref {variance22} since $\mathcal U(1)=Z$ and $\mathcal U(s)=\mathcal F_1(s)$ (these were defined in \eqref {variance}). All we're left to do is prove Lemma \ref {j}.
\begin {center}
\section {Proof of Lemma \ref {j}}
\end {center}
Here we prove Lemma \ref {j}. We suspect the exponent may be pushed down to 1 which would be the best input for the fourth moment and we return to this another time. For now in this article we're content with the simplest saving, the main point being the cancellation between integrals explained in the introduction.
\\
\\ Let $r$ be as in \eqref {defn} so that 
\begin {equation}\label {r}
r=1/p+\mathcal O(1/p^2)
\end {equation}
Let $\Gamma _\Delta (l)$ be as in \eqref {defn} and for $\sigma >1$ and a Dirichlet character $\chi $ let
\begin {eqnarray}\label {idef}
U(\Delta )&:=&\prod _{p|\Delta }\frac {1}{1+\frac {2r}{p}}\hspace {10mm}K(p):=1+\frac {r^2U(p)}{p}\hspace {10mm}K_\Delta (l):=\prod _{p|l\atop {p\nmid \Delta }}K(p)\notag 
\\ \mathcal G_\Delta ^{\chi }(s)&:=&\sum _{l=1}^\infty \frac {\Gamma _\Delta (l)\chi (l)}{l^s}\hspace {10mm}\mathcal K_{d,\Delta }(s):=\sum _{l=1\atop {d|l}}^\infty \frac {K_{\Delta }(l)}{l^s}.
\end {eqnarray}
We have\footnote {$\mathcal G_\Delta ^\chi (s)$ is obviously extendable to $\sigma >-1$ but we don't need it} , writing $Y=\chi (p)p^{-s}$,
\begin {align}\label {ii}
\mathcal G_\Delta ^{\chi }(s)&\hspace {1mm}=\hspace {1mm}
L_\chi (s)L_\chi (s+1)\overbrace {\prod _{p\nmid \Delta }
\left (1+rY\right )\left (1-\frac {Y}{p}\right )\prod _{p|\Delta }\left (1-\frac {Y}{p}\right )
}^{\text { holomorphic and $\ll 1$ for $\sigma >-1/2$ by \eqref {r}}}
\end {align}
and in particular, assuming $(d,\Delta )=1$ and writing $Y=p^{-s}$,
\begin {align}\label {iiprincipal}
\mathcal G_\Delta ^{\chi _0}&\hspace {1mm}=\hspace {1mm}L_{\chi _0}(s)L_{\chi _0}(s+1)\prod _{p\nmid d\Delta }\left (1+rY\right )\left (1-\frac {Y}{p^{}}\right )\prod _{p|\Delta }\left (1-\frac {Y}{p^{}}\right )&&\hspace {1mm}=:\hspace {1mm}L_{\chi _0}(s)L_{\chi _0}(s+1)\mathcal R_{d,\Delta }(s)\notag 
\\ &\hspace {1mm}=\hspace {1mm}L_{\chi _0}(s)\prod _{p\nmid d\Delta }\left (1+rY\right )&&\hspace {1mm}=:\hspace {1mm}L_{\chi _0}(s)\mathcal Q_{d\Delta }(s)
\end {align}
and we have, using $K_{\Delta }(ld)=K_\Delta (d)K_{d\Delta }(l)$,
\begin {eqnarray}\label {kk}
\mathcal K_{d}(s)
&=&\frac {\zeta (s)K_\Delta (d)}{d^s}\prod _{p\nmid d\Delta }\left (1+\underbrace {\frac {r^2U(p)Y}{p}}_{\ll 1/p\text { for }\sigma >-2}\right )
\end {eqnarray}
so
\begin {eqnarray}\label {movingintegrals}
\text {for $\sigma >-1/2$}\hspace {6mm}\mathcal G_\Delta ^{\chi }(s)&\ll &\left |L_\chi (s)L_\chi (s+1)\right |\notag 
\\ &\ll &\left |\zeta (s)\zeta (s+1)\right |\max \left \{ d^{-\sigma },1\right \} \hspace {10mm}\text { if }\chi =\chi _0\notag 
\\ &\ll &\max \left \{ (dt)^{1/2-\sigma },1\right \} \hspace {22mm}\text { always }\notag 
\\ \text {for $\sigma >-2$}\hspace {6mm}\mathcal K_{d}(s)&\ll &\left |\frac {\zeta (s)}{d^s}\right |.
\end {eqnarray}
Let, where $\Gamma $ is the Gamma function, $C_8$ is as in \eqref {defn}, and $\Gamma _\Delta ,\mathcal Q_{d\Delta },\mathcal R_{1,\Delta }$ are as above,
\begin {eqnarray}\label {bc}
A_\Delta (d)&:=&\frac {2\Gamma (3)\mathcal Q_{d\Delta }(1)^2}{\Gamma (6)}\sum _{e|d}\frac {\Gamma _\Delta (e)^2}{\phi (d/e)e^2}L_{\chi _0(d/e)}(1)^2\hspace {10mm}\notag 
\\ c&:=&\frac {2\Gamma (3)\zeta (0)}{\Gamma (5)}\notag 
\\ B_\Delta (d)&:=&\frac {c}{d}\mathcal Q_{d\Delta }(1)\mathcal Q_{d\Delta }(0)\Gamma _\Delta (d)^2\notag 
\\ B_\Delta ^*(d)&:=&\frac {c}{d}\mathcal Q_{d\Delta }(1)\mathcal Q_{d\Delta }(0)\Gamma _\Delta (d)^2\sum _{p|d}\frac {\log p}{\Gamma _\Delta (p)^2}\notag 
\\ Z&:=&\frac {\zeta '(0)}{\zeta (0)}-\frac {\Gamma '(5)}{\Gamma (5)}+1\notag
\\ Z_{\Delta }&:=&\frac {\mathcal R_{1,\Delta }'(0)}{\mathcal R_{1,\Delta }(0)}+Z\notag  
\\ M_{d,\Delta }(X)&:=&A_\Delta (d)X^5+B_\Delta (d)\Big (\log X/d+Z_{d\Delta }\Big )X^4+B_\Delta ^*(d)X^4\notag 
\\ I_{d,\Delta }(X)&:=&\mathcal Q_{d\Delta }(1)\int _{(-1/2)}\frac {\mathcal Q_{d\Delta }(s)}{s(s+2)(s+3)(s+4)}\left (\sum _{e|d}\frac {\Gamma _\Delta (e)^2}{\phi (d/e)e^{s+1}}L_{\chi _0(d/e)}(1)L_{\chi _0(d/e)}(s)\right )X^{s+4}ds\notag 
\\ J_{d,\Delta }(X)&:=&C_8U(\Delta )\int _{(-1/2)}\frac {\mathcal K_{d,\Delta }(s)X^{s+4}ds}{(s+2)(s+3)(s+4)}
\end {eqnarray}
and note in particular, keeping one integral where it is and moving the other to the left,
\begin {eqnarray}\label {bounds}
(d\leq X)\hspace {15mm}I_{d,\Delta }(X)\ll \frac {X^4}{d}\hspace {15mm}J_{d,\Delta }(X)\ll d^2X^{2}.
\end {eqnarray}
Whenever the letters $l,l'$ are in context write $\tilde l=l+l'$. With $r,\Gamma _\Delta $ as in \eqref {defn} define multiplicative $g_\Delta $ through $\Gamma _\Delta =g_\Delta \star 1$ so that
\begin {equation}\label {inverse}
g_\Delta (d)=\left \{ \begin {array}{ll}g_1(d)&\text { if }(d,\Delta )=1\\ 0&\text { if }(d,\Delta )>1\text { or if $d$ is squarefree}\end {array}\right .\hspace {5mm}\text { and }\hspace {5mm}g_1(p)=r
\end {equation}
and so that
\begin {eqnarray}\label {step1}
\sum _{l+l'\leq X}(X-\tilde l)^2\tilde l\Gamma _\Delta (l)\Gamma _\Delta (l')\Gamma _\Delta (\tilde l)&=&\sum _{d\leq X}g_\Delta (d)\sum _{l+l'\leq X\atop {\tilde l\equiv 0(d)}}(X-\tilde l)^2\tilde l\Gamma _\Delta (l)\Gamma _\Delta (l')\notag 
\\ &=:&\sum _{d\leq X}g_\Delta (d)\mathcal L_{d,\Delta }(X).
\end {eqnarray}
Let $C_8,\mathcal U(s)$ be as in \eqref {defn} and let
\begin {eqnarray}\label {vvv0}
V_s(\Delta )&:=&\prod _{p|\Delta }V_s(p):=\prod _{p|\Delta }\left (1+\frac {2r}{p}+\frac {rY}{p^{}}\left (p+3r+r^2\right )\right )\hspace {10mm}Y:=p^{-s}
\end {eqnarray}
so that
\begin {eqnarray}\label {vvv}
\mathcal V_\Delta (s)&:=&\prod _{p\nmid \Delta }V_s(p)\hspace {50mm}\sigma >0\notag 
\\ &=&\frac {\zeta (s+1)}{V_s(\Delta )}\prod _{p}V_s(p)(1-Y/p)\hspace {20mm}\sigma >-1/2\hspace {5mm}(\text {by \eqref {r}})
\end {eqnarray}
so that with \eqref {defn} 
\begin {eqnarray}\label {ayran}
\sum _{\Delta =1}^\infty \frac {I(\Delta )\mathcal V_\Delta (s)}{\Delta ^{s+1}}=\prod _{p}\left (V_s(p)+\frac {YI(p)}{p}\right )=C_8\mathcal U(s)\hspace {40mm}\sigma >0.
\end {eqnarray}
We will prove that
\begin {eqnarray}\label {method1}
\mathcal L_{d,\Delta }(X)&=&M_{d,\Delta }(X)+2\Gamma (3)I_{d,\Delta }(X)+\mathcal O\left (X^\epsilon \left (X^3+d^{3/2}X^{5/2}+d^2X^{2}\right )\right )
\end {eqnarray}
and
\begin {eqnarray}\label {method2}
\mathcal L_{d,\Delta }(X)&=&\frac {C_8U(\Delta )\mathcal K_{d}(1)}{30}X^5+2J_{d,\Delta }(X)+\mathcal O\left (\frac {X^4}{d}\right )
\end {eqnarray}
which from \eqref {bounds} imply
\begin {eqnarray*}
\mathcal L_{d,\Delta }(X)&=&M_{d,\Delta }(X)+2\Gamma (3)I_{d,\Delta }(X)+2J_{d,\Delta }(X)+\mathcal O\left (X^\epsilon \min \left \{ X^3+d^{3/2}X^{5/2}+X^{2}d^{2},\frac {X^4}{d}\right \} \right )
\end {eqnarray*}
and therefore
\begin {eqnarray}\label {gokcen}
\sum _{d=1}^\infty g_\Delta (d)\mathcal L_{d,\Delta }(X)
&=&\sum _{d=1}^\infty g_\Delta (d)M_{d,\Delta }(X)+\sum _{d=1}^\infty g_\Delta (d)\Big (2\Gamma (3)I_{d,\Delta }(X)+2J_{d,\Delta }(X)\Big )+\mathcal O\left (X^{17/5}\right )
\end {eqnarray}
so with \eqref {ayran} we get Lemma \ref {j} if we prove \eqref {method1}, \eqref {method2} and
\begin {eqnarray}\label {series1}
\sum _{d=1}^\infty g_\Delta (d)I_{d,\Delta }(X)&=&\int _{(-1/2)}\frac {\zeta (s)\mathcal V_\Delta (s)X^{s+4}ds}{s(s+2)(s+3)(s+4)}
\end {eqnarray}
\begin {eqnarray}\label {series2}
\sum _{d=1}^\infty g_\Delta (d)J_{d,\Delta }(X)&=&\int _{(0)}\frac {\zeta (s)\mathcal V_\Delta (s)X^{s+4}ds}{(s+2)(s+3)(s+4)}
\end {eqnarray}
\begin {eqnarray}\label {mainterms}
\sum _{\Delta ,d=1}^\infty \Delta ^3I(\Delta )g_\Delta (d)M_{d,\Delta }(X/\Delta )&=&M(X).
\end {eqnarray}
\subsection {Proof of \eqref {method1}}
Recall that the Gamma function is holomorphic for $\sigma >-1/2$ except for a simple pole with residue 1 at $s=0$ and that $L_{\chi _0}(s+1)$ has a simple pole at $s=0$.
\\ 
\\ For a squarefree character of order $d$, we know that $L_{\chi _0}(s)$ has a zero at $s=0$ of order  at least 2 if $d$ has more than one prime factor, and a simple zero if $d=p$ (and no zero if $d=1$), so (recalling the definitions in \eqref {iiprincipal})
\begin {eqnarray}\label {residue}&&Res_{s=0}\left \{ L_{\chi _0}(s)L_{\chi _0}(s+1)\mathcal R_{d,\Delta }(s)\frac {\Gamma (s)+\Gamma (s+1)}{\Gamma (s+5)}X^{s}\right \} \notag 
\\ &&\hspace {10mm}=\hspace {4mm}\underbrace {Res_{s=0}\left \{ \frac {\zeta (s)\zeta (s+1)\mathcal R_{1,\Delta }(s)\Gamma (s)X^s}{\Gamma (s+5)}\right \} +\frac {\zeta (0)\mathcal R_{1,\Delta }(0)}{\Gamma (5)}}_{d=1}\notag 
\\ &&\hspace {14mm}+\hspace {4mm}\underbrace {\frac {\mathcal R_{d,\Delta }(0)}{\Gamma (5)}Res_{s=0}\left \{ L_{\chi _0}(s)L_{\chi _0}(s+1)\Gamma (s)\right \} }_{d=p}\notag 
\\ &&\hspace {10mm}=\hspace {4mm}\underbrace {\frac {\zeta (0)\mathcal Q_{\Delta }(0)}{\Gamma (5)}\left (\frac {\zeta '(0)}{\zeta (0)}+\frac {\mathcal R_{1,\Delta }'(0)}{\mathcal R_{1,\Delta }(0)}-\frac {\Gamma '(5)}{\Gamma (5)}+\log X+1\right )}_{d=1}+\underbrace {\frac {\zeta (0)\log d\hspace {1mm}\mathcal Q_{d\Delta }(0)}{\Gamma (5)}}_{d=p}\notag 
\\ &&\hspace {10mm}=:\hspace {4mm}\underbrace {C^*_\Delta (X)}_{d=1}+\underbrace {C^{**}_{d,\Delta }}_{d=p}.
\end {eqnarray}
Write $\Gamma _\Delta (\mathbf l)=\Gamma _\Delta (l)\Gamma _\Delta (l')$. Splitting the sum according to the value of $e=(l,d)$ and using $\Gamma _{\Delta }(el)=\Gamma _\Delta (e)\Gamma _{\Delta e}(l)$ we have
\begin {eqnarray}\label {dechra2}
\sum _{l+l'\leq X\atop {\tilde l\equiv 0(d)}}(X-\tilde l)^2l\Gamma _\Delta (\mathbf l)
&=&\sum _{e|d}e^3\Gamma _\Delta (e)^2\sum _{l+l'\leq X/e\atop {\tilde l\equiv 0(d/e)\atop {(\mathbf l,d/e)=1}}}(X/e-\tilde l)^2l\Gamma _{\Delta e}(\mathbf l)\notag 
\\ &=&\sum _{e|d}\frac {e^3\Gamma _\Delta (e)^2}{\phi (d/e)}\sum _{\chi (d/e)}\chi (-1)\sum _{l+l'\leq X/e}(X/e-\tilde l)^2l\Gamma _{\Delta e}(\mathbf l)\chi (l)\overline {\chi }(l')\notag 
\\ &=:&\sum _{e|d}\frac {e^3\Gamma _\Delta (e)^2}{\phi (d/e)}\sum _{\chi (d/e)}\chi (-1)\mathcal L_{\Delta e}^{\chi (d/e)}(X/e);
\end {eqnarray}
here and throughout a ``$\chi (q)$" in the superscript is just to remind us of the character's modulus. 
With the usual Perron type representations and recalling \eqref {idef} 
\begin {eqnarray}\label {l}
\mathcal L_\Delta ^{\chi (d)}(X)
=\Gamma (3)\int _{(1)}\int _{(1)}\frac {\mathcal G_\Delta ^{\overline {\chi }(d)}(s)\mathcal G_\Delta ^{\chi (d)}(w)\Gamma (s)\Gamma (w+1)X^{s+w+3}dwds}{\Gamma (s+w+4)}=:\Gamma (3)D_\Delta ^{\chi (d)}(X)
\end {eqnarray}
then, dropping the indices on $\mathcal G$ to avoid clutter and reading $\mathcal G(1)=0$ for $\chi \not =\chi _0$ (recall \eqref {ii}),
\begin {eqnarray}\label {banana}
D_\Delta ^{\chi (d)}(X)&=&\int _{(1)}\mathcal G^{}(s)\Gamma (s)\left (\frac {\mathcal G_{}(1)X^{s+4}}{\Gamma (s+5)}+\int _{(0)}\frac {\mathcal G^{}(w)\Gamma (w+1)X^{s+w+3}dw}{\Gamma (s+w+4)}\right )ds\notag 
\\ &=&\frac {\mathcal G_{}(1)^2X^5}{\Gamma (6)}+\mathcal G_{}(1)\int _{(0)}\frac {\mathcal G^{}(s)\Gamma (s)X^{s+4}ds}{\Gamma (s+5)}\notag 
\\ &&+\hspace {4mm}\int _{(0)}\mathcal G^{}(w)\Gamma (w+1)\left (\frac {\mathcal G_{}(1)X^{w+4}}{\Gamma (w+5)}+\int _{(0)}\frac {\mathcal G^{}(s)\Gamma (s)X^{s+w+3}ds}{\Gamma (s+w+4)}\right )dw\notag 
\\ &=&\frac {\mathcal G_{}(1)^2X^5}{\Gamma (6)}+\mathcal G_{}(1)\int _{(0)}\mathcal G^{}(s)\underbrace {\frac {\Gamma (s)+\Gamma (s+1)}{\Gamma (s+5)}}_{=:g(s)}X^{s+4}ds\notag 
\\ &&+\hspace {4mm}\int _{(0)}\int _{(0)}\frac {\mathcal G^{}(s)\mathcal G^{}(w)\Gamma (s)\Gamma (w+1)X^{s+w+3}dwds}{\Gamma (s+w+4)}.
\end {eqnarray}
Remembering \eqref {iiprincipal} in \eqref {residue} the second term is
\begin {eqnarray}\label {be}
&&\mathcal G_{}(1)\left (Res_{s=0}\left \{ \mathcal G_\Delta^{\chi _0(d)}(s)g(s)X^{s}\right \} +\int _{(-1/2)}\mathcal G_\Delta ^{\chi _0(d)}(s)g(s)X^{s}ds\right )X^4\notag 
\\ &&\hspace {25mm}=:\hspace {4mm}\mathcal G(1)\left (\underbrace {C_\Delta ^*(X)}_{d=1}+\underbrace {C^{**}_{d,\Delta }}_{d=p}+T_\Delta ^{\chi _0(d)}(X)\right )X^4\hspace {10mm}
\end {eqnarray}
whilst when summed over $\chi $ the third term is (not forgetting \eqref {movingintegrals} and reading $\mathcal G(0)=0$ for $\chi \not =\chi _0$)
\begin {eqnarray*}
&\leq &\sum _{\chi }\Bigg |\int _{(0)}\mathcal G^{}(s)\Gamma (s)\left (\frac {\mathcal G_{}(0)X^{s+3}}{\Gamma (s+4)}+\int _{(-1/2)}\frac {\mathcal G^{}(w)\Gamma (w+1)X^{s+w+3}dw}{\Gamma (s+w+4)}\right )ds\Bigg |
\\ &\leq &\left |\mathcal G_{}(0)Res_{s=0}\left \{ \frac {\mathcal G^{}(s)\Gamma (s)X^{s+3}}{\Gamma (s+4)}\right \} \right |+\left |\mathcal G_{}(0)\int _{(-1/2)}\frac {\mathcal G^{}(s)\Gamma (s)X^{s+3}ds}{\Gamma (s+4)}\right |
\\ &&+\hspace {4mm}\sum _\chi \Bigg |\int _{(-1/2)}\mathcal G(w)\Gamma (w+1)\hspace {1mm}Res_{s=0}\left \{ \frac {\mathcal G(s)\Gamma (s)X^{s+w+3}ds}{\Gamma (s+w+4)}\right \} dw
\\ &&\hspace {4mm}+\hspace {4mm}\int _{(-1/2)}\int _{(-1/2)}\frac {\mathcal G^{}(s)\mathcal G^{}(w)\Gamma (s)\Gamma (w+1)X^{s+w+3}dwds}{\Gamma (s+w+4)}\Bigg |
\\ &\ll &\left |L_{\chi _0}(0)^2X^3\right |\log X+\sum _{\chi }\int _{(-1/2)}\left |L_{\chi }(0)L_\chi (1)L_{\chi }(s)L_{\chi }(s+1)\frac {\Gamma (s)}{\Gamma (s+4)}X^{s+3}\right |ds
\\ &&+\hspace {4mm}\sum _\chi \int _{(-1/2)}\int _{(-1/2)}\left |L_\chi (s)L_\chi (s+1)L_\chi (w)L_\chi (w+1)\frac {\Gamma (s)\Gamma (w+1)}{\Gamma (s+w+4)}X^{s+w+3}\right |dwds
\\ &\ll &X^\epsilon \Big (X^3+d^{5/2}X^{5/2}+d^3X^{2}\Big )
\end {eqnarray*}
(here we have comfortably bounded the double integral with a well-known mean value estimate) so from \eqref {step1}, \eqref {dechra2} and \eqref {l} 
\begin {eqnarray*}
\mathcal L_{d,\Delta }(X)&=&2\Gamma (3)\sum _{e|d}\frac {e^3\Gamma _\Delta (e)^2}{\phi (d/e)}\Bigg (\frac {\mathcal G_{\Delta e}^{\chi _0(d/e)}(1)^2(X/e)^5}{\Gamma (6)}+\mathcal G_{\Delta e}^{\chi _0(d/e)}(1)\left (\underbrace {C_{\Delta e}^*(X/e)}_{d/e=1}+\underbrace {C^{**}_{d/e,\Delta e}}_{d/e=p}+T_{\Delta e}^{\chi _0(d/e)}(X/e)\right )(X/e)^4\Bigg )
\\ &&+\hspace {4mm}\mathcal O\left (X^{\epsilon }\sum _{e|d}\frac {e^3\Gamma _\Delta (e)^2}{\phi (d/e)}\Big ((X/e)^3+(d/e)^{3/2}(X/e)^{5/2}+(d/e)^{3}(X/e)^2\Big )\right )
\end {eqnarray*}
and \eqref {method1} follows from the definitions in \eqref {iiprincipal}, \eqref {bc}, \eqref {residue}, \eqref {banana}, \eqref {be} and \eqref {inverse}. 
\subsection {Proof of \eqref {method2}}
Recall the definition of $g_\Delta $ from just before \eqref {inverse} and of $\mathcal L_{d,\Delta }(X)$ from \eqref {step1}. As
\[ \sum _{l+l'=\tilde l\atop {\mathbf D|\mathbf l}}1=\underbrace {\frac {\tilde l}{[D,D']}}_{(D,D')|\tilde l}+\mathcal O(1)\]
we have 
\begin {eqnarray*}
\mathcal L_{d,\Delta }(X)&=&\sum _\mathbf Dg_\Delta (\mathbf D)\sum _{\tilde l\leq X\atop {\tilde l\equiv 0(d)}}(X-\tilde l)^2\tilde l\sum _{l+l'=\tilde l\atop {\mathbf D|\mathbf l}}1
\\ &=&\sum _{\tilde l\leq X\atop {\tilde l\equiv 0(d)}}(X-\tilde l)^2\tilde l^2\sum _{\mathbf D\atop {(D,D')|\tilde l}}\frac {g_\Delta (\mathbf D)}{[D,D']}+\mathcal O\left (\frac {X^4}{d}\right )
\end {eqnarray*}
and then after a straightforward calculation the $\mathbf D$-sum is ($C_8$ is defined in \ref {defn} and $U(\Delta )$ in \eqref {idef})
\begin {eqnarray*}
\prod _p\left (1+\frac {2g_\Delta (p)}{p}\right )\prod _{p|\tilde l}\left (1+\frac {g_\Delta (p)^2}{p\left (1+\frac {2g_\Delta (p)}{p}\right )}\right )
=C_8U(\Delta )K_\Delta (\tilde l)
\end {eqnarray*}
from \eqref {inverse} so, with $\mathcal K_d(s)$ as in \eqref {idef}, a standard Perron representation gives
\begin {eqnarray*}
\mathcal L_{d,\Delta }(X)&=&
2C_8U(\Delta )\int _{(1)}\frac {\mathcal K_{d}(s)X^{s+4}ds}{(s+2)(s+3)(s+4)}+\mathcal O\left (\frac {X^4}{d}\right )\notag 
\\ &=&\frac {C_8U(\Delta )\mathcal K_d(1)X^5}{30}+2C_8U(\Delta )\int _{(-1/2)}\frac {\mathcal K_{d}(s)X^{s+4}ds}{(s+2)(s+3)(s+4)}+\mathcal O\left (\frac {X^4}{d}\right )
\end {eqnarray*}
after \eqref {kk} and this becomes \eqref {method2} after \eqref {bc}.
\subsection {Proofs of \eqref {series1}, \eqref {series2} and \eqref {mainterms}}
With $\Gamma _\Delta $ as in \eqref {defn} define
\begin {equation}\label {qqq}
q_s(n):=\prod _{p|n}q_s(p):=\prod _{p|n}(1+rY)\hspace {10mm}L_{\chi _0(d)}(s)=:\zeta (s)\Delta _s(d)\hspace {10mm}\tau (d):=\sum _{e|d}\frac {\Delta _s(d/e)\Gamma _\Delta (e)^2}{e^{s}}
\end {equation}
so that from \eqref {iiprincipal}
\begin {eqnarray}\label {balwn}
\mathcal Q_{\Delta }(s)&=&\prod _{p\nmid \Delta }q_s(p)\hspace {10mm}\sigma >0\notag 
\\ \mathcal Q_{\Delta }(s)&=&\zeta (s+1)\prod _{p\nmid \Delta }q_s(p)(1-Y/p)\prod _{p|\Delta }(1-Y/p)\hspace {10mm}\sigma >-1\hspace {10mm}Y:=p^{-s}\notag 
\\ \mathcal Q_{d\Delta }(s)&=&
\frac {\mathcal Q_\Delta (s)}{q_s(d)}
\end {eqnarray}
and
\begin {eqnarray}\label {balwn2}
\tau (p)
=1+rY(2+r)\hspace {5mm}\text {for $p\nmid \Delta $}.
\end {eqnarray}
Let $V_s(p),\mathcal V_\Delta (s)$ be as in \eqref {vvv0}, \eqref {vvv} so that \eqref {qqq} gives
\begin {equation}\label {balwn3}
V_0(p)=q_1(p)q_0(p)+\frac {r(1+r)^2}{p}.
\end {equation}
From \eqref {balwn} and \eqref {balwn2}
\begin {eqnarray*}
q_1(p)q_s(p)+\frac {r\tau (p)}{p}=\left (1+\frac {r}{p}\right )\left (1+rY\right )+\frac {r+r^2Y(2+r)}{p}=V_s(p)\hspace {10mm}\text {for $p\nmid \Delta $} 
\end {eqnarray*}
so from \eqref {balwn} and \eqref {inverse}
\begin {eqnarray}\label {sul}
\sum _{d=1}^\infty g_\Delta (d)\mathcal Q_{d\Delta }(1)\mathcal Q_{d\Delta }(s)\sum _{e|d}\frac {\Gamma _\Delta (e)^2}{\phi (d/e)e^{s+1}}L_{\chi _0(d/e)}(1)L_{\chi _0(d/e)}(s)&=&\mathcal Q_{\Delta }(1)\mathcal Q_{\Delta }(s)\zeta (s)\sum _{d=1}^\infty \frac {g_\Delta (d)\tau (d)}{dq_1(d)q_s(d)}\notag 
\\ (\text {initially only for }\sigma >0)\hspace {20mm}&=&\zeta (s)\prod _{p\nmid \Delta }\left (q_1(p)q_s(p)+\frac {r\tau (p)}{p}\right )\notag 
\\ &=&\zeta (s)\mathcal V_\Delta (s)
\end {eqnarray}
so since ($I_{d,\Delta }(X)$ is defined in \eqref {bc})
\begin {eqnarray}\label {t}
\sum _{d=1}^\infty g_\Delta (d)I_{d,\Delta }(X)&=&\int _{(-1/2)}\frac {X^{s+4}}{s(s+2)(s+3)(s+4)}\hspace {1mm}\Big (\text { LHS of \eqref {sul} }\Big )\hspace {1mm}ds
\end {eqnarray}
we get \eqref {series1}. Writing 
\[ l=r+\frac {r^2(1+r)U(p)}{p}\]
we have from \eqref {defn} and \eqref {idef}
\begin {eqnarray*}
K_\Delta (p)\Gamma _\Delta (p)-1&=&\left \{ \begin {array}{ll}l&\text { if }p\nmid \Delta \\ 0&\text { if }p|\Delta \end {array}\right .
\end {eqnarray*}
so from \eqref {idef} and the definition of $g_\Delta $ (this was just before \eqref {inverse}) 
\begin {eqnarray*}
\sum _{d=1}^\infty g_\Delta (d)\mathcal K_{d,\Delta }(s)=\sum _{l=1}^\infty \frac {K_\Delta (l)\Gamma _\Delta (l)}{l^s}
=\zeta (s)\prod _p(1+lY)\prod _{p|\Delta }\frac {1}{1+lY}
\end {eqnarray*}
and (with $V_s(p)$ as in \eqref {vvv0}) 
\begin {eqnarray*}
\frac {1+lY}{U(p)}
=V_s(p)
\end {eqnarray*}
so (with $C_8$ as in \eqref {defn} and $\mathcal V_\Delta (s)$ as in \eqref {vvv})
\begin {eqnarray*}
C_8U(\Delta )\sum _{d=1}^\infty g_\Delta (d)\mathcal K_{d,\Delta }(s)=
\zeta (s)\mathcal V_\Delta (s)
\end {eqnarray*}
so from \eqref {bc} we get
\begin {eqnarray*}
\sum _{d=1}^\infty g_\Delta (d)J_{d,\Delta }(X)&=&\int _{(-1/2)}\frac {\mathcal \zeta (s)\mathcal V_\Delta (s)X^{s+4}ds}{(s+2)(s+3)(s+4)}
\end {eqnarray*}
which is \eqref {series2}. Recall the definitions of $A(d),B_\Delta (d)$ from \eqref {bc}. Putting $s=1$ in \eqref {sul} and using \eqref {ayran} gives 
\begin {eqnarray}\label {bir}
\sum _{\Delta ,d=1}^\infty \frac {I(\Delta )g_\Delta (d)A_\Delta (d)}{\Delta ^2}&=&
\frac {2\Gamma (3)C_8}{\Gamma (6)}\mathcal U(1)
\end {eqnarray}
whilst with \eqref {balwn} and \eqref {inverse}
\begin {eqnarray}\label {bbc}
\sum _{d=1}^\infty \frac {g_\Delta (d)\mathcal Q_{d\Delta }(1)\mathcal Q_{d\Delta }(0)}{d}\Gamma _\Delta (d)^2e^{-s\log d}=\prod _{p|\Delta }(1-1/p)\prod _{p\nmid \Delta }(1-1/p)\left (q_1(p)q_0(p)+\frac {r(1+r)^2e^{-s\log p}}{p}\right )
\end {eqnarray}
so that in particular we get from \eqref {balwn3} (and definitions \eqref {bc} and \eqref {vvv})
\begin {eqnarray}\label {iki2}
\sum _{d=1}^\infty \frac {g_\Delta (d)\mathcal Q_{d\Delta }(1)\mathcal Q_{d\Delta }(0)}{d}\Gamma _\Delta (d)^2
=\mathcal V_\Delta (0).
\end {eqnarray}
so from \eqref {ayran}
\begin {eqnarray}\label {iki}
\sum _{\Delta ,d=1}^\infty \frac {I(\Delta )g_\Delta (d)B_\Delta (d)}{\Delta }&=&cC_8\mathcal U(0)
\end {eqnarray}
Now define 
\begin {align*}
\theta _1(p)&:=\frac {r}{1+r}-1+
\frac {1}{(1+r)^2}
&\theta _2(p)&:=\frac {1}{(1+r)^2}
\\ A_1(n)&:=\sum _{p|n}\theta _1(p)\log p&A_2(n)&:=\sum _{p|n}\theta _2(p)\log p
\\ M_1(n)&:=\frac {1}{nq_1(n)q_0(n)}\sum _{d\Delta =n}I(\Delta )g_\Delta (d)\Gamma _\Delta (d)^2&M_2(\Delta )&:=\frac {I(\Delta )}{\Delta V_0(\Delta )}
\\ M_1^*(n)&:=c\mathcal Q_1(1)\mathcal Q_1(0)M_1(n)&M_2^*(\Delta )&:=c\mathcal V_\Delta (0)M_2(\Delta )
\\ \mathcal A_1(s)&:=\sum _{n=1}^\infty M_1^*(n)e^{sA_1(n)}&\mathcal A_2(s)&:=\sum _{n=1}^\infty M_2^*(n)e^{sA_1(n)}
\end {align*}
and recall the values of $I,\Gamma _\Delta ,g_\Delta ,q_s$ from \eqref {defn}, \eqref {inverse}, \eqref {qqq}. We have for $p\nmid \Delta $
\begin {eqnarray*}
\frac {M_1(p)}{1+M_1(p)}-\frac {M_2(p)}{1+M_2(p)}&=&\underbrace {\frac {I(p)+g_1(p)\Gamma _1(p)^2}{p(1+r)(1+r/p)+I(p)+g_1(p)\Gamma _1(p)^2}}_{=:\frac {G}{Q+G}}-\frac {I}{Q+G}
\end {eqnarray*}
and therefore, writing $R=1/(1+r)$,
\begin {eqnarray*}
\frac {M_1(p)\theta _1(p)}{1+M_1(p)}-\frac {M_2(p)\theta _2(p)}{1+M_2(p)}&=&-\frac {RG}{Q+G}+\frac {R^2(G-I)}{Q+G}=\frac {r(1+2r)}{(1+r)(p+pr+r)}
\end {eqnarray*}
so
\begin {eqnarray*}
\frac {1}{p-1}-\frac {r}{1+r}+\frac {M_1(p)\theta _1(p)}{1+M_1(p)}-\frac {M_2(p)\theta _2(p)}{1+M_2(p)}&=&\frac {p(1-pr+3r)}{(p-1)(p+pr+r)}
=\frac {1}{p(p-1)}
\end {eqnarray*}
from \eqref {defn}. Also $\mathcal A_2(0)=\mathcal A_1(0)=cC_8\mathcal U(0)$ from \eqref {iki2} and \eqref {ayran}. From \eqref {iiprincipal} 
\begin {eqnarray*}
\frac {\mathcal R_{1,d\Delta }'(0)}{\mathcal R_{1,d\Delta }(0)}&=&
\sum _{p}\left (\frac {1}{p-1}-\frac {r}{1+r}\right )\log p+\sum _{p|d\Delta }\frac {r}{1+r}\log p=:\mathcal S+\sum _{p|d\Delta }\frac {r}{1+r}\log p\hspace {10mm}\text {for }(d,\Delta )=1
\end {eqnarray*}
so from \eqref {bc} (and not forgetting that \eqref {defn} and \eqref {inverse} cause the $d\Delta ,n$ variables below to be all squarefree)
\begin {eqnarray*}
&&\sum _{\Delta ,d=1}^\infty \frac {I(\Delta )g_\Delta (d)}{\Delta }\Bigg (-B_\Delta (d)\log (d\Delta )+B_\Delta (d)Z_{d\Delta }+B^*_\Delta (d)\Bigg )\notag 
\\ &&\hspace {10mm}=\hspace {4mm}\sum _{\Delta ,d=1}^\infty \frac {I(\Delta )g_\Delta (d)\mathcal Q_{d\Delta }(1)\mathcal Q_{d\Delta }(0)}{d\Delta }\Gamma _\Delta (d)^2\notag 
\\ &&\hspace {20mm}\times \hspace {4mm}\left (Z+\mathcal S+\sum _{p|d\Delta }\frac {r\log p}{1+r}-\underbrace {\log (d\Delta )}_{=\sum _{p|d\Delta }\log p}+\sum _{p|d\Delta }\frac {\log p}{(1+r)^2}-\sum _{p|\Delta }\frac {\log p}{(1+r)^2}\right )\notag 
\\ &&\hspace {10mm}=\hspace {4mm}(Z+\mathcal S)\sum _{n=1}^\infty M_1^*(n)+\sum _{n=1}^\infty M_1^*(n)A_1(n)-\mathcal V_1(0)\sum _{\Delta =1}^\infty M_2^*(\Delta )A_2(\Delta )\hspace {15mm}\text {(from \eqref {balwn})}\notag 
\\ &&\hspace {10mm}=\hspace {4mm}(Z+\mathcal S)\mathcal A_1(0)+\mathcal A_1'(0)-\mathcal A_2'(0)\notag 
\\ &&\hspace {10mm}=\hspace {4mm}cC_8\mathcal U(0)\left (Z+\sum _p\Bigg (\frac {1}{p-1}-\frac {r}{1+r}+\frac {M_1(p)\theta _1(p)}{1+M_1(p)}-\frac {M_2(p)\theta _2(p)}{1+M_2(p)}\Bigg \} \log p\right )\notag 
\\ &&\hspace {10mm}=\hspace {4mm}cC_8\mathcal U(0)\left (Z+\sum _p\frac {\log p}{p(p-1)}\right ).
\end {eqnarray*}
which with \eqref {bir} and \eqref {iki} is \eqref {mainterms}.
\begin {center}

\end {center}
\hspace {1mm}
\\
\\
\\
\\
\\  
\\ \emph {Tomos Parry
\\ Bilkent University, Ankara, Turkey
\\ tomos.parry1729@hotmail.co.uk}

\end {document}